\documentclass[12pt,reqno]{amsart}

\usepackage{amssymb}

\allowdisplaybreaks

\numberwithin{equation}{section}

\newtheorem{Theorem}{Theorem}

\newtheorem{Lemma}[Theorem]{Lemma}
\newtheorem{Corollary}[Theorem]{Corollary}

\newtheorem{Conjecture}[Theorem]{Conjecture}

\theoremstyle{definition}

\theoremstyle{remark}

\def\al{\alpha}
\def\be{\beta}

\def\de{\delta}
\def\ep{\varepsilon}

\def\la{\lambda}

\def\bfalpha{\boldsymbol{\alpha}}
\def\bfbeta{\boldsymbol{\beta}}

\begin{document}

\newbox\Adr
\setbox\Adr\vbox{
\centerline{\sc J. M. Brunat$^{1*}$,
C.~Krattenthaler$^{2\dagger\ddagger}$, A.~Lascoux$^{3\dagger}$ and 
A. Montes$^{1*}$}
\vskip18pt
\centerline{$^1$Departament de Matem\`atica Aplicada II,}
\centerline{Universitat Polit\`ecnica de Catalunya,}
\centerline{Jordi Girona 1--3,   08034 Barcelona,   Spain}
\centerline{WWW: \footnotesize\tt http://www-ma2.upc.edu/\~{}montes}
\vskip18pt
\centerline{$^2$Institut Camille Jordan, Universit\'e Claude Bernard
Lyon-I,}
\centerline{21, avenue Claude Bernard, F-69622 Villeurbanne Cedex,
France}
\centerline{WWW: \footnotesize\tt http://igd.univ-lyon1.fr/\~{}kratt}
\vskip18pt
\centerline{$^3$Institut Gaspard Monge, Universit\'e de
Marne-la-Vall\'ee,} 
\centerline{F-77454 Marne-la-Vall\'ee Cedex~2,
France}
\centerline{WWW: \footnotesize\tt http://www-igm.univ-mlv.fr/\~{}al}
}

\title{Some composition determinants}

\author[J. M. Brunat, C.~Krattenthaler, A.~Lascoux and A. Montes]{\box\Adr}

\address{Departament de Matem\`atica Aplicada II,
Universitat Polit\`ecnica de Catalunya,
Jordi Girona 1--3,   08034 Barcelona,   Spain.}

\address{Institut Camille Jordan, Universit\'e Claude Bernard Lyon-I, 
21, avenue Claude Bernard, F-69622 Villeurbanne Cedex, France.}

\address{Institut Gaspard Monge, Universit\'e de Marne-la-Vall\'ee,
F-77454 Marne-la-Vall\'ee Cedex~2, France.}

\thanks{$^*$Work partially supported by the Ministerio de Ciencia y
Tecnolog\'\i a under projects BFM2003-00368 and MTM2004-01728 and
by the Generalitat de Catalunya under project 2001 SGR 00224.\newline
 \indent$^\dagger$Research partially supported by EC's IHRP Programme,
grant HPRN-CT-2001-00272, ``Algebraic Combinatorics in Europe."\newline
 \indent$^\ddagger$Current address: Fakult\"at f\"ur Mathematik, Universit\"at Wien,
 Nordbergstra{\ss}e~15, A-1090 Vienna, Austria}


\begin{abstract}
We compute several parametric determinants in which rows and columns are indexed by
compositions, where the entries are either products of
binomial coefficients or products of
powers. These results generalize previous determinant evaluations due
to the first and fourth author [{\it SIAM J. Matrix Anal.\ Appl.} {\bf 23}
(2001), 459--471] and [``A polynomial generalization of the 
power-compositions determinant," {\it Linear Multilinear Algebra} 
(to appear)], and they prove two conjectures of the second author
[``Advanced determinant calculus: a complement," preliminary version].
\end{abstract}

\maketitle

\section{Introduction}

A composition of a non-negative integer $n$ 
is a vector $(\al_1,\al_2,\dots,\al_k)$ of non-negative integers such that
$\al_1+\al_2+\dots+\al_k=n$, for some $k$. For a fixed $k$, 
let $\mathcal C(n,k)$ denote the corresponding set of compositions of $n$.
While working on a problem in global optimisation, 
two of the authors \cite{BrMoAA}
discovered the following surprising determinant evaluation.
It allowed them to show how to explicitly express a multivariable
polynomial as a {\it difference of convex functions}.
In the statement, we use standard multi-index notation: 
if $\boldsymbol\al=(\al_1,\al_2,\dots,\al_k)$ and
$\boldsymbol \be=(\be_1,\be_2,\dots,\be_k)$ are two compositions, we
let
$$\boldsymbol\al^{\boldsymbol \be}:=\al_1^{\be_1}\al_2^{\be_2}\cdots
\al_k^{\be_k},$$
where $0^0$ is interpreted as $1$.

\begin{Theorem} \label{thm:Brunat}
For any positive integers $n$ and $k$, we have
\begin{equation} \label{eq:Brunat} 
\det_{\boldsymbol\al,\boldsymbol\be\in\mathcal
C(n,k)}\left(\boldsymbol\al^{\boldsymbol \be}\right)=
n^{{\binom{ n + k-1} k}+k-1}\,
 \prod_{i = 1}^{ n-1}
      i^{( n -i+1) {\binom {n + k-i-1} { k-2}}}.
\end{equation}
\end{Theorem}

In the preliminary version \cite{KratBZa} of \cite{KratBZ}, 
the second author observed empirically that there seemed to be
a polynomial generalisation of this theorem.

\begin{Conjecture}[{\cite[Conjecture~57]{KratBZa}}] \label{conj:Brunat1}
For any positive integers $n$ and $k$, we have
\begin{equation} \label{eq:Brunat1} 
\det_{\boldsymbol\al,\boldsymbol\be\in\mathcal
C(n,k)}\left((x+\boldsymbol\al)^{\boldsymbol \be}\right)=
(kx+n)^{{\binom{ n + k-1} k}}n^{k-1}\,
 \prod_{i = 1}^{ n-1}
      i^{( n -i+1) {\binom {n + k-i-1} { k-2}}},
\end{equation}
where $x$ is a variable, and $x+\boldsymbol\al$ is short for
$(x+\al_1,x+\al_2,\dots,x+\al_k)$.
\end{Conjecture}

Morover, he also worked out a binomial version of this conjecture.
Extending our
multi-index notation, let
$$\binom{\boldsymbol\al}{\boldsymbol \be}:=
\binom{\al_1}{\be_1}\binom{\al_2}{\be_2}\cdots
\binom{\al_k}{\be_k}.$$

\begin{Conjecture}[{\cite[Conjecture~58]{KratBZa}}] \label{conj:Brunat2}
For any positive integers $n$ and $k$, we have
\begin{equation} \label{eq:Brunat2a} 
\det_{\boldsymbol\al,\boldsymbol\be\in\mathcal
C(n,k)}\left(\binom{x+\boldsymbol\al+\boldsymbol\be}{\boldsymbol \be}\right)=
 \frac{\displaystyle\prod_{i = 0}^{n-1}{( kx+n+k+i) }^
       {{\binom{ k+i-1} { k-1}}}}{\displaystyle\prod_{i = 1}^{n}
      i^{{\binom{n+k-i-1} {k-1}}}},
\end{equation}
where $x$ is a variable, and $x+\boldsymbol\al+\boldsymbol\be$ 
is short for
$(x+\al_1+\be_1,x+\al_2+\be_2,\dots,x+\al_k+\be_k)$.
\end{Conjecture}

In the recent paper \cite{BrMoAB}, the first and fourth author 
succeeded to prove
Conjecture~\ref{conj:Brunat1}. In fact, they established 
the following multivariable generalisation.

\begin{Theorem} \label{thm:Brunat3}
Let $\mathbf x=(x_1,x_2,\dots,x_k)$ be a vector of indeterminates.
Then, for any positive integers $n$ and $k$, we have
\begin{equation} \label{eq:Brunat3} 
\det_{\boldsymbol\al,\boldsymbol\be\in\mathcal
C(n,k)}\left((\mathbf x+\boldsymbol\al)^
    {\boldsymbol \be}\right)=
 \left({ \vert \mathbf x\vert+ n }\right)
      ^{\binom{n+k-1} {k}}
\prod _{i=1} ^{n}i^{(k-1)\binom {n+k-i-1}{k-1}},
\end{equation}
where $\mathbf x+\boldsymbol\al$ is short for
$(x_1+\al_1,x_2+\al_2,\dots,x_k+\al_k)$, and
where $\vert \mathbf x\vert=x_1+x_2+\dots+x_k$. 
\end{Theorem}

The purpose of this paper is to, in some sense, explain the
miraculous existence of all these formulae. We do this by introducing
further variables, $\la_1,\la_2,\dots,\la_k$, 
in the binomial determinant in
\eqref{eq:Brunat2a}, and by proving an evaluation theorem for the
resulting determinant. All the afore-mentioned determinant evaluations
are then special cases, respectively limit cases, of this new theorem.
To be precise, the main result of this paper is the following
determinant evaluation.

\begin{Theorem} \label{thm:Brunat4}
Let $\mathbf x=(x_1,x_2,\dots,x_k)$ and
$\boldsymbol\la=(\la_1,\la_2,\dots,\la_k)$ be vectors of
indeterminates.
Then, for any positive integers $n$ and $k$, we have
\begin{equation} \label{eq:Brunat4} 
\det_{\boldsymbol\al,\boldsymbol\be\in\mathcal
C(n,k)}\left(\binom{\mathbf x+\boldsymbol\la\boldsymbol\al}
    {\boldsymbol \be}\right)=\left(\prod _{j=1}
^{k}\la_j\right)^{\binom {n+k-1}k}
\frac {\displaystyle\underset{\ep_1+\dots+\ep_k<n}
{\prod _{\ep_1,\dots,\ep_k\ge0} ^{}} \left(
n+ \sum _{j=1} ^{k}\left(\frac {x_j} {\la_j}-\frac {\ep_j}
{\la_j}\right) \right)}
{\displaystyle\prod _{i=1} ^{n}i^{\binom{n+k-i-1} {k-1}}},
\end{equation}
where $\mathbf x+\boldsymbol\la\boldsymbol\al$ is short for
$(x_1+\la_1\al_1,x_2+\la_2\al_2,\dots,x_k+\la_k\al_k)$.
\end{Theorem}

Using the elementary property $\binom Xm=(-1)^m\binom {-X+m-1}m$
of binomial coefficients, we show
that an equivalent way to write the same result is as follows.

\begin{Theorem} \label{thm:Brunat4b}
Let $\mathbf x=(x_1,x_2,\dots,x_k)$ and
$\boldsymbol\la=(\la_1,\la_2,\dots,\la_k)$ be vectors of
indeterminates.
Then, with notation as in Theorem~{\em\ref{thm:Brunat4}}, 
for any positive integers $n$ and $k$ we have
\begin{equation} \label{eq:Brunat4b} 
\det_{\boldsymbol\al,\boldsymbol\be\in\mathcal
C(n,k)}\left(\binom{\mathbf x+\boldsymbol\la\boldsymbol\al+\boldsymbol\be}
    {\boldsymbol \be}\right)=\left(\prod _{j=1}
^{k}\la_j\right)^{\binom {n+k-1}k}
\frac {\displaystyle\underset{\ep_1+\dots+\ep_k<n}
{\prod _{\ep_1,\dots,\ep_k\ge0} ^{}} \left(n+
 \sum _{j=1} ^{k}\left(\frac {x_j+1} {\la_j}+\frac {\ep_j}
{\la_j}\right) \right)}
{\displaystyle\prod _{i=1} ^{n}i^{\binom{n+k-i-1} {k-1}}}.
\end{equation}
\end{Theorem}

In order to see how Conjecture~\ref{conj:Brunat2} is implied by these
results, it is convenient to first state separately the special cases of
Theorems~\ref{thm:Brunat4} and \ref{thm:Brunat4b} where all the
$\la_i$'s are identical.
In that case, the product in the numerator on the right-hand sides of
\eqref{eq:Brunat4} and \eqref{eq:Brunat4b} can be rearranged by
grouping together the factors corresponding to compositions
$\boldsymbol \varepsilon=(\ep_1,\ep_2,\dots,\ep_k)$ of $n-i$, which become
identical. Taking into account that the number of such compositions
is $\binom {n-i+k-1}{k-1}$, this yields the following two
corollaries.

\begin{Corollary} \label{thm:Brunat2}
Let $\mathbf x=(x_1,x_2,\dots,x_k)$ be a vector of indeterminates,
and let $\la$ be an indeterminate.
Then, for any positive integers $n$ and $k$, we have
\begin{equation} \label{eq:Brunat2} 
\det_{\boldsymbol\al,\boldsymbol\be\in\mathcal
C(n,k)}\left(\binom{\mathbf x+\la\boldsymbol\al}
    {\boldsymbol \be}\right)=\la^{(k-1)\binom {n+k-1}k}
\prod _{i=1} ^{n} \left(\frac{{ \vert \mathbf x\vert+(\la-1)n+i }}
 i\right)
      ^{\binom{n+k-i-1} {k-1}},
\end{equation}
where $\mathbf x+\la\boldsymbol\al$ is short for
$(x_1+\la\al_1,x_2+\la\al_2,\dots,x_k+\la\al_k)$, and
where $\vert \mathbf x\vert=x_1+x_2+\dots+x_k$, as before. 
\end{Corollary}

\begin{Corollary} \label{thm:Brunat2b}
Let $\mathbf x=(x_1,x_2,\dots,x_k)$ be a vector of indeterminates,
and let $\la$ be an indeterminate.
Then, with notation as in Corollary~{\em\ref{thm:Brunat2}}, 
for any positive integers $n$ and $k$ we have
\begin{equation} \label{eq:Brunat2b} 
\det_{\boldsymbol\al,\boldsymbol\be\in\mathcal
C(n,k)}\left(\binom{\mathbf x+\la\boldsymbol\al+\boldsymbol\be}
    {\boldsymbol \be}\right)=\la^{(k-1)\binom {n+k-1}k}
\prod _{i=1} ^{n} \left(\frac{{ \vert \mathbf x\vert+(\la +1)n+k-i }}
 i\right)
      ^{\binom{n+k-i-1} {k-1}}.
\end{equation}
\end{Corollary}

Clearly, Conjecture~\ref{conj:Brunat2} is the special case of the above
corollary where $\la=1$ and $x_i=x$ for all $i$. 

Theorem~\ref{thm:Brunat3} is also implied by
Theorem~\ref{thm:Brunat4}. To see this, we shall show that, by 
extracting the highest homogeneous component in \eqref{eq:Brunat4}
(this could also be realised by an appropriate limit),
we obtain the following corollary.

\begin{Corollary} \label{thm:Brunat5a}
Let $\mathbf x=(x_1,x_2,\dots,x_k)$ and
$\boldsymbol\la=(\la_1,\la_2,\dots,\la_k)$ be vectors of
indeterminates.
Then, with notation as in Theorem~{\em\ref{thm:Brunat4}}, 
for any positive integers $n$ and $k$ we have
\begin{equation} \label{eq:Brunat5a} 
\det_{\boldsymbol\al,\boldsymbol\be\in\mathcal
C(n,k)}\left((\mathbf x+\boldsymbol\la\boldsymbol\al)^
    {\boldsymbol \be}\right)=\left(\prod _{j=1}
^{k}\la_j\right)^{\binom {n+k-1}k}
\left(n+ \sum _{j=1} ^{k}\frac {x_j} {\la_j}\right)^{\binom {n+k-1}k}
\prod _{i=1} ^{n}i^{(k-1)\binom {n+k-i-1}{k-1}}.
\end{equation}
\end{Corollary}

The corresponding special case where all the $\la_i$'s are identical
is the following.

\begin{Corollary} \label{thm:Brunat1a}
Let $\mathbf x=(x_1,x_2,\dots,x_k)$ be a vector of indeterminates,
and let $\la$ be an indeterminate.
Then, with notation as in Corollary~{\em\ref{thm:Brunat2}}, 
for any positive integers $n$ and $k$ we have
\begin{equation} \label{eq:Brunat1a} 
\det_{\boldsymbol\al,\boldsymbol\be\in\mathcal
C(n,k)}\left((\mathbf x+\la\boldsymbol\al)^
    {\boldsymbol \be}\right)=\la^{(k-1)\binom {n+k-1}k}
 \left({ \vert \mathbf x\vert+\la n }\right)
      ^{\binom{n+k-1} {k}}
\prod _{i=1} ^{n}i^{(k-1)\binom {n+k-i-1}{k-1}}.
\end{equation}
\end{Corollary}

Clearly, Theorem~\ref{thm:Brunat3} is the special case $\la=1$ of this
corollary.

In the next section, we give proofs of Theorems~\ref{thm:Brunat4}
and \ref{thm:Brunat4b}, and of
Corollary~\ref{thm:Brunat5a} (and, thus, of
Corollaries~\ref{thm:Brunat2}, \ref{thm:Brunat2b} and
\ref{thm:Brunat1a} also). 
%We provide in fact two proofs
%of our main theorem, Theorem~\ref{thm:Brunat4}. 
%Our first proof is based on the ``identification of factors" technique
%(see \cite[Sec.~2.4]{KratBN}). In contrast, our second proof
%extends the inductive procedure from
%\cite{BrMoAA,BrMoAB} that was used in the original proofs of the special cases
%in Theorems~\ref{thm:Brunat} and \ref{thm:Brunat3}.
In contrast to the inductive procedure in
\cite{BrMoAA,BrMoAB} that was used in the original proofs of Theorems~\ref{thm:Brunat} and \ref{thm:Brunat3},
our proof is based on the ``identification of factors" technique
(see \cite[Sec.~2.4]{KratBN}).
As it turns out, the crucial identity
in both of our proofs is the multivariate version of the Chu--Vandermonde
summation formula (see Lemma~\ref{lem:ChuVand}).
Finally, in the last section, we derive analogues of
Theorems~\ref{thm:Brunat4} and \ref{thm:Brunat4b}, and
of Corollary~\ref{thm:Brunat5a} for the
subdeterminants in which we restrict the rows and columns to 
compositions of $n$ with exactly
$k$ {\it positive} summands.

\section{The proofs}

\begin{Lemma} \label{lem:ChuVand}
Let $\mathbf x=(x_1,x_2,\dots,x_k)$ be a vector of indeterminates, and
let $n$ and $k$ be non-negative integers. Then
$$\sum _{\boldsymbol\de\in\mathcal C(n,k)} ^{}
\binom {\mathbf x}{\boldsymbol\de} =
\sum _{\de_1+\dots+\de_k=n} ^{}\binom {\mathbf x}{\boldsymbol\de} =
\binom {\vert \mathbf x\vert}n.$$
\end{Lemma}

\begin{proof}
The Chu--Vandermonde summation formula (see e.g.\ \cite[Sec.~5.1,
(5.27)]{GrKPAA}) reads
$$
\sum _{r=0} ^{s}\binom Mr\binom {N}{s-r}=\binom {M+N}s.$$
On the basis of this formula, the assertion of the lemma is
easily proved by induction on $k$.
\end{proof}

\begin{proof}%[First proof of Theorem~\ref{thm:Brunat4}]
[Proof of Theorem~\ref{thm:Brunat4}]
We prove the theorem by the identification of factors method
described in \cite[Sec.~2.4]{KratBN}.
For convenience, let us write $M(n,k)$ for the matrix of which we want
to compute the determinant, that is,
$$M(n,k)=\left(\binom{\mathbf x+\boldsymbol\la\boldsymbol\al}
    {\boldsymbol \be}\right)_{\boldsymbol\al,\boldsymbol\be\in\mathcal
C(n,k)}.$$

\medskip
{\it Step~1. The term
$$T(n,k,\mathbf x,\boldsymbol\la,\boldsymbol\ep):=n\prod _{j=1} ^{k}\la_j+
\sum _{t=1} ^{k} {x_t}\underset{j\ne t}{\prod _{j=1} ^{k}}\la_j
-\sum _{t=1} ^{k}  {\ep_t} \underset{j\ne t}{\prod _{j=1} ^{k}}\la_j$$
divides $\det M(n,k)$ for any composition $\boldsymbol\ep\in\mathcal
C(n-i,k)$, $1\le i\le n$.} (It should be noted that
$T(n,k,\mathbf x,\boldsymbol\la,\boldsymbol\ep)$ is the factor
corresponding to $\ep_1,\ep_2,\dots,\ep_k$ in the product 
in the numerator on the
right-hand side of \eqref{eq:Brunat4}, up to multiplication by $
\prod _{j=1} ^{k}\la_j$.) To prove this assertion, we
find a vector in the kernel of
\begin{equation} \label{eq:Mnk}
M(n,k)\Big\vert_{T(n,k,\mathbf x,\boldsymbol\la,\boldsymbol\ep)=0}.
\end{equation}
This kernel lives in the free vector space generated by the
compositions in $\mathcal C(n,k)$. Given  $\boldsymbol \de\in\mathcal
C(n,k)$, let us denote the corresponding
element (``unit vector") in this vector space by $e_{\boldsymbol\de}$.
Then we claim that the vector
\begin{equation} \label{eq:v} 
v_{\boldsymbol\ep}:=
\sum _{\boldsymbol\de\in\mathcal C(i,k)} ^{}
\left(\sum _{t=1} ^{k}\frac {\de_t} {\la_t}\right)
\binom {\boldsymbol \de+\boldsymbol \ep}{\boldsymbol \de}
e_{\boldsymbol\de+\boldsymbol\ep}
\end{equation}
is in the kernel of the matrix \eqref{eq:Mnk}. To
see this, we calculate, using Lemma~\ref{lem:ChuVand} and the
notation $\mathbf u_t=(0,\dots,0,1,0,\dots,0)$ (with the $1$ in
position $t$),
\begin{align*}
\text{coefficient of $e_\al$ in }M(n,k)\cdot v_{\boldsymbol\ep}&=\sum
_{\boldsymbol\de\in\mathcal C(i,k)} ^{} 
\binom{\mathbf x+\boldsymbol\la\boldsymbol\al}
    {\boldsymbol \de+\boldsymbol \ep}
\left(\sum _{t=1} ^{k}\frac {\de_t} {\la_t}\right)
\binom{\boldsymbol\de+\boldsymbol\ep}
    {\boldsymbol \de}\\
&\kern-3cm
=\binom {\mathbf x+\boldsymbol\la\boldsymbol\al}
{\boldsymbol \ep}
\sum _{t=1} ^{k}\sum
_{\vert\boldsymbol\de\vert=i} ^{} \frac {\de_t} {\la_t}
\binom{\mathbf x+\boldsymbol\la\boldsymbol\al-\boldsymbol \ep}
    {\boldsymbol \de}\\
&\kern-3cm
=\binom {\mathbf x+\boldsymbol\la\boldsymbol\al}
{\boldsymbol \ep}
\sum _{t=1} ^{k}\frac {1} {\la_t}(x_t+\la_t\al_t-\ep_t)
\sum
_{\vert\boldsymbol\de\vert=i} ^{} 
\binom{\mathbf x+\boldsymbol\la\boldsymbol\al-\boldsymbol \ep-\mathbf u_t}
    {\boldsymbol \de-\mathbf u_t}\\
&\kern-3cm
=\binom {\mathbf x+\boldsymbol\la\boldsymbol\al}
{\boldsymbol \ep}
\left(n+\sum _{t=1} ^{k}\left(\frac {x_t} {\la_t}-
\frac {\ep_t} {\la_t}\right)\right)
\binom{\vert\mathbf x+\boldsymbol\la\boldsymbol\al\vert-(n-i)-1}
    {i-1}.
\end{align*}
Since $i\ge1$, the occurrence of the factor in the middle implies
$$M(n,k)\Big\vert_{T(n,k,\mathbf x,\boldsymbol\la,\boldsymbol\ep)=0}\cdot
v_{\boldsymbol\ep}=0.$$

\medskip
{\it Step~2. Comparison of degrees.}
By inspection,
the (total) degree in the $x_i$'s and $\la_i$'s 
of the determinant on the left-hand side of
\eqref{eq:Brunat4} is at most $n\cdot\vert \mathcal C(n,k)\vert=n\binom
{n+k-1}n$.
%\begin{align*}
%\sum _{\boldsymbol\be\in\mathcal C(n,k)} ^{}\be_1&=
%\sum _{i=0} ^{n}i\cdot\vert\mathcal C(n-i,k-1)\vert\\&=
%\sum _{i=0} ^{n}i\binom {n-i+k-2}{k-2}\\
%&=-\sum _{i=0} ^{n}(n-i)\binom {n-i+k-2}{k-2}+
%\sum _{i=0} ^{n}n\binom {n-i+k-2}{k-2}\\
%&=-(k-1)\sum _{i=0} ^{n}\binom {n-i+k-2}{k-1}+
%n\sum _{i=0} ^{n}\binom {n-i+k-2}{k-2}\\
%&=-(k-1)\binom {n+k-1}{k}+n\binom {n+k-1}{k-1}\\
%&=\binom {n+k-1}{k}.
%\end{align*}
%Here, we used special instances of the Chu--Vandermonde summation to
%evaluate the sums over $i$.
On the other hand, the degree in the $x_i$'s and $\la_i$'s 
of the right-hand side of
\eqref{eq:Brunat4} is equal to $k\binom {n+k-1}k=n\binom {n+k-1}n$.
%\begin{equation} \label{eq:vand1} 
%\sum _{i=1} ^{n}\binom {n+k-i-1}{k-1}=\binom {n+k-1}k.
%\end{equation}
%Again, we used a special instance of the Chu--Vandermonde summation to
%evaluate the sum over $i$.

Since the degree bound on the determinant is the same as the degree
of the right-hand side of \eqref{eq:Brunat4}, 
the determinant must be equal to the
right-hand side up to a possible multiplicative constant
which is independent of the $x_i$'s and $\la_i$'s. 
%Moreover, since, by symmetry, the same
%is also true for $x_i$, $i=2,3,\dots,k$, 
We conclude that
\vbox{
\begin{multline} \label{eq:const}
\det_{\boldsymbol\al,\boldsymbol\be\in\mathcal
C(n,k)}\left(\binom{\mathbf x+\boldsymbol\la\boldsymbol\al}
    {\boldsymbol \be}\right)=\text{const}(n,k)
\left(\prod _{j=1}
^{k}\la_j\right)^{\binom {n+k-1}k}\\
\times
{\displaystyle\underset{\ep_1+\dots+\ep_k<n}
{\prod _{\ep_1,\dots,\ep_k\ge0} ^{}} \left(
n+ \sum _{j=1} ^{k}\left(\frac {x_j} {\la_j}-\frac {\ep_j}
{\la_j}\right) \right)},
\end{multline}}

\noindent
where $\text{const}(n,k)$ is independent of $\mathbf x$ and $\boldsymbol\la$. 

\medskip
{\it Step~3. Computation of the multiplicative constant.}
If set $x_i=0$ and $\la_i=1$ for all $i$, then the determinant on the
left-hand side of \eqref{eq:Brunat4} becomes triangular with $1$s on
the diagonal. Thus, we obtain that
$$1=\text{const}(n,k)\underset{\ep_1+\dots+\ep_k<n}
{\prod _{\ep_1,\dots,\ep_k\ge0} ^{}} \left(
n - \sum _{j=1} ^{k}{\ep_j}\right) =
\text{const}(n,k)\prod _{i=1} ^{n} i^{\vert\mathcal C(n-i,k)\vert}.
$$
This implies that
$$\text{const}(n,k)=\prod _{i=1} ^{n} i^{-\binom {n+k-i-1}{k-1}},$$
completing the proof of the theorem.
\end{proof}

\begin{proof}[Proof of the equivalence of Theorem~\ref{thm:Brunat4}
and \ref{thm:Brunat4b}]
If we replace $\la_i$ by $-\la_i$ and $x_i$ by $-x_i-1$ for all $i$ in
\eqref{eq:Brunat4}, and then use the
identity $\binom Xm=(-1)^m\binom {-X+m-1}m$, then we obtain
\begin{align}\notag 
&\det_{\boldsymbol\al,\boldsymbol\be\in\mathcal
C(n,k)}\left((-1)^n
\binom{\mathbf x+\la\boldsymbol\al+\boldsymbol\be}
    {\boldsymbol \be}\right)=\\
\label{eq:-1}
&
=(-1)^{k\binom {n+k-1}k}
\left(\prod _{j=1}
^{k}\la_j\right)^{\binom {n+k-1}k}
\frac {\displaystyle\underset{\ep_1+\dots+\ep_k<n}
{\prod _{\ep_1,\dots,\ep_k\ge0} ^{}} \left(
n+ \sum _{j=1} ^{k}\left(\frac {x_j+1} {\la_j}+\frac {\ep_j}
{\la_j}\right) \right)}
{\displaystyle\prod _{i=1} ^{n}i^{\binom{n+k-i-1} {k-1}}}.
\end{align}
Except for the signs, this is exactly \eqref{eq:Brunat4b}. However,
we have
\begin{multline*}
 n\cdot\vert \mathcal C(n,k)\vert+k\binom {n+k-1} k
= n\binom {n+k-1} { k-1}+k\binom {n+k-1} k
=2\frac{(n-k-1)!}{(k-1)!\,(n-1)!}.
\end{multline*}
Since this is an even number, the signs in \eqref{eq:-1} do indeed cancel.
\end{proof}

\begin{proof}[Proof of Corollary~\ref{thm:Brunat5a}]
The right-hand and left-hand sides of \eqref{eq:Brunat4} are both
polynomials in the $\la_i$'s and the $x_i$'s. As we already observed in
Step~2 of the first proof of Theorem~\ref{thm:Brunat4},
the (total) degree in the $\la_i$'s and the
$x_i$'s of the determinant on the left-hand side is equal to
$$n\cdot\vert\mathcal C(n,k)\vert=n\binom {n+k-1}{k-1}=
k\binom {n+k-1}k,$$ 
and the degree in the $\la_i$'s and the
$x_i$'s of the expression on the right-hand side is exactly the same
value.
Therefore, if we extract the homogeneous parts in the $\la_i$'s and the
$x_i$'s of degree $k\binom {n+k-1}k$ in \eqref{eq:Brunat4}, we obtain
\begin{equation*}
\det_{\boldsymbol\al,\boldsymbol\be\in\mathcal
C(n,k)}\left(\frac {(\mathbf x+\boldsymbol\la\boldsymbol\al)^{\boldsymbol\be}}
    {\boldsymbol \be!}\right)=\left(\prod _{j=1}
^{k}\la_j\right)^{\binom {n+k-1}k}
\frac {\displaystyle\underset{\ep_1+\dots+\ep_k<n}
{\prod _{\ep_1,\dots,\ep_k\ge0} ^{}} \left(
n+ \sum _{j=1} ^{k}\frac {x_j} {\la_j} \right)}
{\displaystyle\prod _{i=1} ^{n}i^{\binom{n+k-i-1} {k-1}}},
\end{equation*}
where $\boldsymbol \be!=
\prod _{i=1} ^{k}\be_i!$, or, equivalently,
\begin{equation*}
\det_{\boldsymbol\al,\boldsymbol\be\in\mathcal
C(n,k)}\left( {(\mathbf x+\la\boldsymbol\al)^{\boldsymbol\be}}
\right)=\left(\prod _{j=1}
^{k}\la_j\right)^{\binom {n+k-1}k}
\frac{\left(n+ \sum _{j=1} ^{k}\frac {x_j} {\la_j}\right)^{\sum _{i=0}
^{n-1}\vert\mathcal C(i,k)\vert}}
{\prod _{i=1} ^{n}i^{\binom{n+k-i-1} {k-1}}}
\Bigg(
\prod _{\boldsymbol\be\in\mathcal
C(n,k)} ^{}\boldsymbol \be!\Bigg).
\end{equation*}
Since 
$$\sum _{i=0} ^{n-1}\vert\mathcal C(i,k)\vert=
\sum _{i=0} ^{n-1}\binom {i+k-1}{k-1}=\binom {n+k-1}k,$$
the only missing piece for the proof of the corollary is the
verification of the identity
\begin{equation} \label{eq:beta}
 \prod _{\boldsymbol\be\in\mathcal
C(n,k)} ^{}\boldsymbol \be!=
\prod _{\be_1+\dots+\be_k=n} ^{}\be_1!\cdots\be_k!=
\prod _{i=1} ^{n} i
      ^{k\binom{n+k-i-1} {k-1}}.
\end{equation}
This can, for example, be done by induction on $n+k$, by using the
obvious recurrence
$$\Pi(n,k)=
\prod _{i=0} ^{n}\left(i!^{\vert\mathcal C(n-i,k-1)}\Pi(n-i,k-1)\right),$$
where
$$\Pi(n,k)=\prod _{\be_1+\dots+\be_k=n} ^{}\be_1!\cdots\be_k!.$$
\end{proof}

\section{Determinants for compositions with only positive parts}

Let $C^*(n,k)$ denote the set of all compositions of $n$ with
exactly $k$ positive summands.
Then we have the following theorem.

\begin{Theorem}
\label{DP}
Let $\mathbf x=(x_1,x_2,\dots,x_k)$ and
$\boldsymbol\la=(\la_1,\la_2,\dots,\la_k)$ be vectors of
indeterminates.
Then, for any positive integers $n$ and $k$, $n\ge k$, we have
\begin{multline} \label{eq:DP}
\det_{\bfalpha,\bfbeta\in C^*(n,k)}
\left(\binom{\mathbf{x}+\boldsymbol\lambda\bfalpha}{\bfbeta} \right)=
\left(\prod _{j=1}
^{k}\la_j\right)^{\binom {n-1}k}
\left(
\prod_{i=1}^{k}\prod_{j=1}^{n-k+1}
\left(\frac { x_i+\lambda_i j}j\right)^{\binom {n-j-1} { k-2}}
\right)\\
\times
\frac {\displaystyle\underset{\ep_1+\dots+\ep_k<n-k}
{\prod _{\ep_1,\dots,\ep_k\ge0} ^{}} \left(
n+ \sum _{j=1} ^{k}\left(\frac {x_j-1} {\la_j}-\frac {\ep_j}
{\la_j}\right) \right)}
{\displaystyle\prod _{i=1} ^{n}i^{\binom{n-i-1} {k-1}}}.
\end{multline}
\end{Theorem}
\begin{proof}
Let $\mathbf{1}$ be the $k$-vector with all entries equal to 1.
The mapping $C^*(n,k)\rightarrow C(n-k,k)$ defined by 
$\boldsymbol{\alpha}=(\alpha_1,\ldots,\alpha_k)
\mapsto \boldsymbol{\alpha}-\mathbf{1}
=(\alpha_1-1,\ldots, \alpha_k-1)$
is bijective. Thus, we have
\begin{align*}
\det_{\boldsymbol{\alpha},\boldsymbol{\beta}\in C^*(n,k)}
&\left(\binom {\mathbf{x}+\boldsymbol\lambda\boldsymbol{\alpha}}
{\boldsymbol{\beta}}\right)  = 
\det_{\boldsymbol{\alpha},\boldsymbol{\beta}\in C^*(n,k)}
\left( 
\binom {\mathbf{x}+\boldsymbol\lambda+\boldsymbol\lambda(\boldsymbol{\alpha}-\mathbf{1})}
{\boldsymbol{\beta}-\mathbf{1}+\mathbf{1}}
\right) \\
&= 
\det_{\boldsymbol{\alpha},\boldsymbol{\beta}\in C(n-k,k)}
\left(
\binom {\mathbf{x}+\boldsymbol\lambda+\boldsymbol\lambda\boldsymbol{\alpha}}
{\boldsymbol{\beta}+\mathbf{1}}
\right)  \\
&=  
\det_{\boldsymbol{\alpha},\boldsymbol{\beta}\in C(n-k,k)}
\left(
 \binom{\mathbf{x}+\boldsymbol\lambda-\mathbf{1}+\boldsymbol\lambda\bfalpha}{\bfbeta}
\prod _{i=1} ^{k}\frac {
x_i+\lambda_i+\lambda_i\alpha_i}{\beta_i+1}
\right)  \\
&=
\left(
\prod_{\boldsymbol{\alpha}\in C(n-k,k)}
\prod _{i=1} ^{k}\frac {
x_i+\lambda_i(\alpha_i+1)}{\alpha_i+1}
\right)\det_{\boldsymbol{\alpha},\boldsymbol{\beta}\in C(n-k,k)}
\left(
 \binom{\mathbf{x}+\boldsymbol\lambda-\mathbf{1}+\boldsymbol\lambda\bfalpha}{\bfbeta}
\right).
\end{align*}
In the product, a ratio $(x_i+\lambda_i (j+1))/(j+1)$ 
appears as many times as there are
compositions in $\mathcal C(n-k-j,k-1)$, i.e., ${\binom {n-j-2} { k-2}}$ times.
Thus,
\begin{align*}
\prod_{\boldsymbol{\alpha}\in C(n-k,k)}
\prod _{i=1} ^{k}\frac {
x_i+\lambda_i(\alpha_i+1)}{\alpha_i+1}
&=\prod_{i=1}^k\prod_{j=0}^{n-k} 
\left(\frac {x_i+\lambda_i(j+1)} {j+1}\right)^{\binom {n-j-2} { k-2}}\\
&=
\prod_{i=1}^k\prod_{j=1}^{n-k+1} \left(\frac 
{x_i+\lambda_i j}j\right)^{\binom {n-j-1} { k-2}}.
\end{align*}
Substituting this and using Theorem~\ref{thm:Brunat4} with $n$ replaced 
by $n-k$ and $\mathbf x$ replaced by $\mathbf x+\boldsymbol\la-\mathbf1$ 
to evaluate the determinant, we obtain the desired formula. 
\end{proof}  

By replacing $\la_i$ by $-\la_i$ and $x_i$ by $-x_i-1$ for all $i$ in
\eqref{eq:DP}, and then using the
identity $\binom Xm=(-1)^m\binom {-X+m-1}m$, we obtain the following
equivalent form of Theorem~\ref{DP}.

\begin{Theorem}
\label{DP2}
Let $\mathbf x=(x_1,x_2,\dots,x_k)$ and
$\boldsymbol\la=(\la_1,\la_2,\dots,\la_k)$ be vectors of
indeterminates.
Then, for any positive integers $n$ and $k$, $n\ge k$, we have
\begin{multline} \label{eq:DP2}
\det_{\bfalpha,\bfbeta\in C^*(n,k)}
\left(\binom{\mathbf{x}+\boldsymbol\lambda\bfalpha+\bfbeta}{\bfbeta} \right)=
\left(\prod _{j=1}
^{k}\la_j\right)^{\binom {n-1}k}
\left(
\prod_{i=1}^{k}\prod_{j=1}^{n-k+1}
\left(\frac { x_i+\lambda_i j+1}j\right)^{\binom {n-j-1} { k-2}}
\right)\\
\times
\frac {\displaystyle\underset{\ep_1+\dots+\ep_k<n-k}
{\prod _{\ep_1,\dots,\ep_k\ge0} ^{}} ^{} \left(
n+ \sum _{j=1} ^{k}\left(\frac {x_j+2} {\la_j}+\frac {\ep_j}
{\la_j}\right) \right)}
{\displaystyle\prod _{i=1} ^{n}i^{\binom{n-i-1} {k-1}}}.
\end{multline}
\end{Theorem}

Extracting the highest homogeneous component in Theorem~\ref{DP}, we obtain
the following analogue of Corollary~\ref{thm:Brunat1a}.

\begin{Corollary} \label{thm:Brunat1b}
Let $\mathbf x=(x_1,x_2,\dots,x_k)$ and
$\boldsymbol\la=(\la_1,\la_2,\dots,\la_k)$ be vectors of
indeterminates.
Then, for any positive integers $n$ and $k$, $n\ge k$, we have
\begin{multline}
\det_{\bfalpha,\bfbeta\in C^*(n,k)}
\left((\mathbf{x}+\boldsymbol\lambda\bfalpha)^{\bfbeta} \right)=
\left(\prod _{j=1}
^{k}\la_j\right)^{\binom {n-1}k}
\left(
\prod_{i=1}^{k}\prod_{j=1}^{n-k+1}
\left( { x_i+\lambda_i j}\right)^{\binom {n-j-1} { k-2}}
\right)\\
\times
\left(
{\displaystyle\underset{\ep_1+\dots+\ep_k<n-k}
{\prod _{\ep_1,\dots,\ep_k\ge0} ^{}} \left(
n+ \sum _{j=1} ^{k}\frac {x_j} {\la_j}\right)}\right)
\prod _{i=1} ^{n-k}i^{(k-1)\binom {n-i-1}{k-1}}.
\end{multline}
\end{Corollary}

\end{document}